%% file: AMA.tex
\documentclass{amsart}
\usepackage{enumerate}
\usepackage{amssymb,amsmath}
\usepackage{latexsym}
\usepackage{graphicx}
\usepackage[]{color}
\usepackage{ifthen}
\usepackage{tcolorbox}

\newcommand{\cA}{\mathcal{A}}
\newcommand{\fA}{\mathfrak{A}}
\newcommand{\cF}{\mathcal{F}}

\newcommand{\cL}{\mathcal{L}}

\newcommand{\R}{\mathbb{R}}
\newcommand{\Forb}{\mathrm{Forb}}

\newtheorem{theorem}{Theorem}[section]

\newtheorem{proposition}[theorem]{Proposition}
\newtheorem{corollary}[theorem]{Corollary}
\newtheorem{problem}{Problem}[section]
\newtheorem{conjecture}{Conjecture}[]

\theoremstyle{definition}
\newtheorem{definition}[theorem]{Definition}
\newtheorem{example}[theorem]{Example}
\newtheorem{examples}[theorem]{Examples}

\theoremstyle{remark}

\numberwithin{equation}{section}
\newif\ifmargin
\margintrue
\newif\ifskip
\skiptrue

\newcommand{\ZZ}{\mathbb{Z}}
\newcommand{\NN}{\mathbb{N}}
\newcommand{\PP}{\mathbb{P}}
\newcommand{\EE}{\mathbb{E}}
\newcommand{\CC}{\mathbb{C}}
\newcommand{\cC}{\mathcal{C}}

\newcommand{\aA}{\mathcal{A}}

\newcommand{\bP}{\mathbf{P}}
\newcommand{\bNP}{\mathbf{NP}}
\newcommand{\goesto}{\rightarrow}

\newcommand{\SOL}{\mathrm{SOL}}
\newcommand{\MSOL}{\mathrm{MSOL}}
\newcommand{\FPT}{\mathrm{FPT}}
\newcommand{\DC}{\mathrm{DC}}
\newcommand{\DCE}{\mathrm{DCE}}
\newcommand{\TG}{\mathrm{TG}}
\newcommand{\VEC}{\mathrm{VEC}}
\newcommand{\RR}{\mathrm{RR}}
\newcommand{\CMSOL}{\mathrm{CMSOL}}

\begin{document}

\title{Meta-theorems for Graph Polynomials}

\author{Johann A. Makowsky}
\address{Department of Computer Science, Technion - IIT, Haifa, Israel}
\email{janos@cs.technion.ac.il}
\thanks{Based on my talks given at the MATRIX-Institute workshop, September 2023}



\date{May 10, 2024}


\begin{abstract}
In this paper I survey the sources of inspiration for my own and co-authored work 
in trying to develop a general theory of graph polynomials.
I concentrate on meta-theorems, i.e., theorem which depend only on the form infinite classes of graph polynomials are
defined in some formalism, and not on the actual meaning of the particular definitions. 
\end{abstract}

\maketitle
\small
\tableofcontents
\normalsize
\input{AMA-short}

\input{AMA-prereq}
\input{NMA-recurrence}

\input{NMA-complexity}
\input{MA-unimodality}

\input{MA-harary}
\input{MA-weakly}

\input{NMA-semantic}

\input{NMA-conclu}

\bibliographystyle{alpha}
\bibliography{Quotes,MA}
\end{document}

%% file: AMA-short.tex
\section{Introduction}

\subsection{One, two, many graph polynomials}
In my paper \cite{makowsky2008zoo}
I set myself the goal of developping a general theory of graph polynomials.
The paper was written after my realization, together with Boris Zilber, that there are many more naturally occurring
graph polynomials
than are treated in the literature, see my papers \cite{makowsky2023got} and  
Kotek, Makowsky and Zilber \cite{ar:KotekMakowskyZilber08}. 
In analogy to the matching polynomial or the independence polynomial, which use the number of $k$-matchings or 
$k$-independent sets as coefficients of a generating function, 
one can use
the number of configurations on $k$ vertices (or edges) given by a description in some logical formalism
as coefficients of a generating function. 
The domination polynomial is a good example, but there are no limits in creating such graph polynomials.
We also noticed that counting the number of conditional colorings treated in the literature
gives graph polynomials in the number of colors which I named later Harary polynomials, see Section \ref{se:harary}.
Earlier attempts to define large classes of graph polynomials go back to Farrell \cite{ar:Farrell79a}.

The abundance of naturally defined graph polynomials led some mathematicians in search of a topic to 
define somehow an arbitrarily defined new graph polynomial, say $F(G;x)$, and investigate it
by dealing with one or more of the following questions:
\begin{enumerate}[(i)]
\item
Compute $F(G_n;x)$ for selected graph families $G_n: n \in \NN$;
\item
Find $F(G;x)$-unique graphs;
\item
Interpret $F(G;a)$ for selected values $a \in \ZZ$ in terms of graph theoretic properties of $G$; 
\item
Interpret the polynomials $F(G_n;x)$ for certain sequences of graphs $G_n$;
\item
Compare the distinguishing power of $F(G;x)$ to the distinguishing power of some other graph polynomial $F'(G;x)$.
\item
Is $F(G;x)$ in some sense ``equivalent'' to some previously studied graph polynomial $F_0(G;x)$?
\item
Where are the roots of $F(G;x)$ located?
\end{enumerate}
The same questions can also be asked for multivariate polynomials.

Some of the more exciting papers in graph polynomials have dealt with such questions, but more recently
I have seen many papers where the choice of $F(G;x)$ is badly, if at all, motivated, and only a few answers to
these questions were given.

\subsection{Meta-theorems}
The purpose of this paper is to show what I had, and still have, in mind 
by a ``general theory of graph polynomialsa'' 
in \cite{makowsky2008zoo}. 

Let me first give a (admittedly trivial) example of a meta-theorem.

Assume we want to show that a certain graph property, say triangle-free, is hereditary (closed under induced subgraphs).§
We can show this directly in the framework of graphs. But we can also observe that triangle-free can be expressed
by a universal first order formula $\phi$:
$$
\phi :
\forall u,v ( E(u,v) \leftrightarrow E(v,u)) \wedge
\forall u,v,w ( 
\neg E(u,v) \vee
\neg E(v,w)) \vee
\neg E(w,u))
).
$$
Now it is easy to show that any graph property expressible by a universal first order formula $\psi$
is hereditary. 
This is a {\em meta-theorem} since it uses only the fact that $\psi$ is a universal first order formula,
irrespective of the particular meaning of $\psi$.
If we allow infinite graphs as well, the well known {\L}o\'s-Tarski  Theorem asserts that the 
converse is true as well: 
\begin{theorem}[{\L}o\'s-Tarski Theorem]
Let $\cA$ be a class of
(possible infinite) graphs axiomtizable by a set of first order sentences $\Sigma_{\cA}$.
Then $\cA$ is hereditary if and only if $\Sigma_{\cA}$ is logically equivalent to a set of universal first order sentences.
\end{theorem}
If restricted to finite structures, Tait \cite{Tait1959} showed that it is wrong even if $\Sigma_{\cA}$ is finite.

There is another way of formulating a meta-theorem for hereditary classes of graphs.
A class of finite graphs  $\cA$ is hereditary if and only if it can be characterized by a (possibly infinite) set 
$\Forb_ {\cA}$ of forbidden induced subgraphs. There the mere existence of $\Forb_ {\cA}$, irrespective
of the elements of $\Forb_ {\cA}$ guarantees that $\cA$ is hereditary.
This can also be formulated in infinitary logic by a universal formula of $\mathcal{L}_{\infty, \omega}$, see
\cite{rosen1995preservation}.
Applied to graph polynomials rather than graph properties
we should look for theorems of the following type:
\begin{tcolorbox}
Let $\cF$ be a class of graph polynomials 
and $\cC$ a class of graphs, both
defined in some logical formalism. 
Then for every $G \in \cC$ and $F(G;\bar{x}) \in \cF$
the graph polynomial $F(G;\bar{x})$ exhibits a characteristic behaviour or property
irrespective of the  particular meaning of the definition of $\cF$ and $\cC$.
\end{tcolorbox}
This can be viewed  as a meta-theorem if it only depends on the form of how $\cF$
and $\cC$ are defined, and not on the particular choice of these definitions.

An early, possibly the earliest, meta-theorem in graph theory is the Specker-Blatter Theorem from 1981, 
\cite{blatter1981lenombre,specker1990application,fischer2011application}.

\begin{theorem}[Specker-Blatter Theorem]
Let $\phi$ is a formula in Monadic Second Order Logic in the language of graphs.
Let $S_{\phi}(n)$ be number of binary relation $E \subset [n]^2$ such that
graph $G = ([n], E)$ satisfies $\phi$.
Then sequence
$S_{\phi}^m(n) = S_{\phi}(n) \pmod{m}$
is ultmately periodic with
$$
S_{\phi}^m(n) = S_{\phi}^m(n+p_{m,\phi})
$$
for all $n \geq q_{m,\phi}$.
\end{theorem}
We consider this a meta-theorem, since ultimate periodicity is a consequence of the syntactic form of $\phi$ only
(but the values $p_{m,\phi}, q_{m,\phi}$ depend on $m$ and $\phi$).
At the time of publication in 1981 (and till 2000) this theorem was widely unnoticed.
For its history and recent developments, see \cite{fischer2024extensions}.

\subsection{Towards a general theory of graph polynomials}
Since  the publication of \cite{makowsky2008zoo} my collaborators and I have achieved some progress in developing the 
modest beginnings of such a general theory.
In the following sections we give examples and discuss such meta-theorems in detail.

In Section \ref{se:recursion} we discuss theorems inspired by results about the Tutte polynomial
and Tutte-Grothendieck invariants ($\TG$-invariants) and similar theorems.
More generally, we outline our results from \cite{godlin2012graph} about unwinding recursive definition of
graph polynomials.

In Section \ref{se:complexity} we review general theorems about the complexity of evaluating
graph polynomials. This includes Fixed Parameter Tractable ($\FPT$) cases for graph classes of bounded
width (tree-width, clique-width, rank-width, etc), and dichotomy theorems for graph polynomials
arising from partition functions.

In Section \ref{se:unimodality} we discuss general theorems about real-rooted graph polynomials
and graph polynomials where the coefficients satisfy concavity conditions.
The special cases of graph polynomials where the coefficients are unimodal is here of special interest.

In Section \ref{se:harary} we review some results about Harary polynomials.

In Section \ref{se:weakly} we discuss a conjecture by Bollob\'as, Pebody and Riordan which states that
the Tutte polynomial is almost complete. A graph polynomial $F(G;\bar{x})$ is almost complete
if almost all graphs are $F$-unique.
We contrast this to graph polynomials $F(G;\bar{x})$ being weakly distinguishing if almost
all graphs $G$ have a $F$-mate and prove a meta-theorem about weakly distinguishing graph polynomials.

In Section \ref{se:semantic} We briefly introduce and discuss the distinction between syntactic and semantic statements 
about (families of) graph polynomials.

Finally, in Section \ref{se:conclu} we draw conclusions and suggest some further research along the lines described
in this paper.

\small
\subsection*{Acknowledgements}
My PhD students of the last 25 years, Udi Rotics \cite{phd:Rotics}, Ilia Averbouch \cite{phd:Averbouch}, 
Tomer Kotek \cite{phd:Kotek} and Vsevolod Rakita \cite{phd:Rakita}, contributed
much to these developments. So did my co-authors, Benny Godlin, Orli Herscovici, Emilia Katz, Elena Ravve 
and Ruixue Zhang.
Markus Bl\"aser, 
Holger Dell, 
Andrew Goodall,
Miki Hermann and
Steve Noble joined me in studying the complexity of graph polynomials.
I would like to thank Graham Farr, Kerri Morgan, Vsevolod Rakita, Elena Ravve and Peter Tittmann for reading and commenting
on earlier drafts of this paper. Their observations and sugeestions were all incorporated into the final version.

\subsection*{The MATRIX Institute}
Two events of the MATRIX Institute had great impact on my work described in this paper:
the Tutte Centenary Retreat in December 2017 and the Workshop on Uniqueness and Discernement in Graph Polynomials
in October 2023. 
I would like to thank the MATRIX Institute and the organizers of these meetings for their fruitful initiatives:
Graham Farr, Dillon Mayhew, Kerri Morgan, James Oxley and Gordon Royle in 2017, and
Jo Ellis-Monaghan, Iain Moffatt, Kerri Morgan and Graham Farr in 2023.

\normalsize


%% file: AMA-prereq.tex
\section{Prerequisites}
\label{se:prereq}

This paper is not self-contained. We assume the reader is somehow familiar with the concepts listed below.
\begin{description}
\item[\bf Graphs]
We consider a  simple graph as relational structure $G = (V(G), E(G))$ where $V(G)$ is the universe (the set of vertices)
and $E(G) \subseteq V(G)^2$ is a binary relation.
We denote by $n(G) = |V(G)|$ the number of vertices, by $m(G) = |E(G)|$ the number of edges,
and $k(G)$ the number of connected components of $G$.
$G$ may have loops, directed edges, but no multiple edges.
Let $e \in E(G)$ be an edge of $G$. We denote by $G_{-e}$, $G_{/e}$ and $G_{\dagger e}$
the graph obtained from $G$ by, respectively, removing, contracting  and extracting the edge $e$.
Edge extraction removes not only $e = (u,v)$ but also the vertices $u, v \in V(G)$.

If we allow multiple edges, we consider the graph as an hypergraph $H =(V(H), E(H), inc)$
with two disjoint sets of vertices $V(H)$ and edges $E(H)$ and a binary incidence relation $inc \subseteq V(H) \times E(H)$
with $inc(u,e)$ if and only if $u$ is a vertex of $e$.

A graph property is {\em hereditary} if its closed under induced subgraphs,
{\em monotone} if its closed under (not necessarily induced) subgraphs, and it is {\em minor-closed}
if it is closed under minors.

\item[\bf Logics]
We assume the reader is familiar with the various fragments of Second Order Logic $\SOL$ such as
Monadic Second Order Logic $\MSOL$, possibly augmented by modular counting quantifiers $\CMSOL$, etc.
For references see \cite{ebbinghaus1999finite,bk:Libkin2004},
In $\SOL$ graphs $G = (V(G), E(G))$ and graphs as hypergraphs $H =(V(H), E(H), inc)$ are bi-interpretable.
In $\MSOL$ or $\CMSOL$ this is not the case. 
We write $\MSOL_g$ for the former and $\MSOL_h$ for the latter (and the same for $\CMSOL$).
This distinction is crucial in Section \ref{se:complexity}  and also in the monograph of Courcelle and Engelfriet 
\cite{courcelle2012graph}.

\item[\bf $\SOL$-definable graph polynomials]
Given a fragment $\cL$ of $\SOL$ we denote by $GP(\cL)$ the family of graph polynomials
definable in $\cL$. The exact formalism is described in \cite{ar:KotekMakowskyZilber11}.
A special univariate case looks like this
$$
F(G;x) = \sum_{A \subset V(G): (G(V), E(G), A) \models \phi} x^{|A|}
$$
or
$$
F(G;x) = \sum_{B \subset E(G): (G(V), E(G), B) \models \phi} x^{|A|}
$$
where $\phi$ is a formula in $\cL$.

\item[\bf Harary polynomials]
Harary  polynomials are of the form
\begin{gather}
\chi_{\cA}(G; x) = \sum_{i \geq 1} b_i^{\cA}(G) x_{(i)},
\end{gather}
where $b_i^{\cA}(G)$ is the number of partitions of $V(G)$ into $i$ non-empty parts, where each part induces a graph in ${\cA}$,
and $x_{(i)}$ is the falling factorial.

Harary polynomials are special cases of generalized chromatic polynomials as described in
\cite{ar:MakowskyZilber2006,ar:KotekMakowskyZilber08,ar:KotekMakowskyZilber11}.
They are the counting analogues of conditional colourings as defined in \cite{harary1985conditional}.

One could look at the case where ${\cA}$ is definable in some fragment of $\SOL$.
However, Theorem \ref{th:main-1} shows that under very conditions on ${\cA}$ the resulting polynomial is not
$\MSOL$-definable.

\item[\bf Complexity classes]
For complexity classes see \cite{bk:papadimitriou94,bk:AroraBarak2009}.
The complexity zoo is a webpage: {\tt https://complexityzoo.net/Complexity\_Zoo},
with a printable version \cite{aaronson2005complexity}.

In this paper in Section \ref{se:complexity} three complexity classes are used: Polynomial time $\bP$, 
the counting class $\sharp\bP$ (sharp $\bP$),
and the complexity class $\FPT$
of Fixed Parameter Tractable problems. 
A parametrized graph problem consists of pairs $A =(G, p(G))$ with $G$ a graph and $p(G)$ a graph parameter of $G$.
We say that $A$ is in $\FPT$ for $p(G)$ if membership in $A$  can be solved in
time $f(p(G)) n^d$ for a function $f(p(G))$ 
and a constant $d \in \NN$
not dependent on $G$.

Good references for $\FPT$ are the monographs
\cite{bk:DowneyFellows99,bk:FlumGrohe2006}.
For complexity of graph polynomials
see also \cite{pr:Makowsky06,goodall2018complexity}.

\item[\bf Tree- and other widths]
We assume the reader is familiar with the notions of tree-width and clique-width. 
A survey of various notions of graph width can be found in \cite{hlinveny2008width}.
\end{description}

%% file: NMA-recurrence.tex
\section{Graph polynomials defined by recurrence relations}
\label{se:recursion}
Let $G =(V(G),E(G))$ be a graph. For an edge $e \in E(G)$
we denote by $n(G), m(G), k(G)$ the number of vertices, edges and connected components of a graph $G = (V(G), E(G))$.
For an edge $e \in E(G)$ we denote by $G_{-e}, G_{/e}, G_{\dagger e}$ the graph obtained from $G$
be deleting, contracting or extracting the edge $e$.

\subsection{Universal graph polynomials}
Here we discuss meta-theorems for graph polynomials of the following form:

\begin{tcolorbox}
Let $\cF$ be a class of graph polynomials defined by a recursion. 
Then there exists a universal graph polynomial
$U_{\cF}(G: \bar{x})$ such that every graph polynomial $F(G;\bar{y})$ can be written as
$$
F(G;\bar{y}) = g(n(G), m(G), k(G)) \cdot U_{\cF}(G: \bar{t}(\bar{y})).
$$
In words, $F(G;\bar{y})$ is the product of a function which depends only on $n(G), m(G), k(G)$,
and a substitution instance of $U_{\cF}(G: \bar{x})$.
\end{tcolorbox}
We consider such theorems as meta-theorems when the existence of the universal polynomial
for a class of graph polynomials $\cF$ 
depends only on the form of the recursion in which $\cF$ is defined.
The logical formalism for the meta-theorem is therefore given by the recursions.

A good source for the classical theorems of this form is \cite{aigner2007course}.
The newer results were obtained in oral collaboration with Bruno Courcelle
\cite{ar:Courcelle08}, with my PhD student Ilia Averbouch
\cite{phd:Averbouch}, and with my seminar students Benny Godlin and Emilia Katz \cite{godlin2012graph}.

\subsection{Substitution instances and prefactors}
Two graphs are {\em similar} if the have the same number of {\bf V}ertices $n(G)$, {\bf E}dges $m(G)$ and 
{\bf C}onnected components $k(G)$.
If any of these three graph parameters  differ on two graphs, they are in an obvious way non-isomorphic.
In principle other notions of similarity with respect to a finite set of graph parameters can be considered.
So far this possibility has not yet been  investigated  by us and  in the literature.
but in this paper we disregard this possibility.

A graph polynomial $f(G; \bar{x})$ is a $\VEC$-invariant  if it invariant under similar graphs.
If a graph polynomial $F(G;\bar{x})$ satisfies
$$
F(G;\bar{x}) = f(G;\bar{x}) \cdot F'(G;\bar{x})
$$ 
for some $\VEC$-invariant $f(G;\bar{x})$, we call $f(G;\bar{x})$ a {\em prefactor}.

Let $\bar{x} = (x_1, \ldots, x_k)$.
$F'(G;\bar{x})$ is a {\em substitution instance of $F(G;\bar{x})$} if there are polynomials
$s_i \in \R[\bar{x}]: i \in [k]$ such that
$$
F'(G;\bar{x}) = F(G; s_1(\bar{x}), \ldots, s_k(\bar{x})).
$$
The substitutions can be rational function in the indeterimates,
but it is important that the substitutions do not depend on $G$.

\begin{example}
The {\em bivariate Potts model} is the graph polynomial
$$
Z(G;y,q) = \sum_{A \subseteq V(G)} y^{|A|} q^{k(G[A])}.
$$ 
$G[A]$ is the spanning subgraph of $G$ with edges in $A$.
For the chromatic polynomial we have
$$
\chi(G;x) = Z(G;x, -1).
$$
There is no prefactor $\chi(G;x)$ is a substitution instance of $Z(G;y,q)$.

\end{example}

\begin{example}
Let $T(G;x,y)$ be the Tutte polynomial and $\chi(G;x)$ be the chromatic polynomial.
We have
$$
\chi(G;x) = (-1)^{n(G)-k(G))} x^{k(G)} \cdot T(G; 1-x,0).
$$
Here $f(G;x) = (-1)^{n(G)-k(G))} x^{k(G)} $ is a prefactor and
$T(G; 1-x,0)$ is a substitution instance of $T(G;x,y)$.
\end{example}

\subsection{Universality}
Let $\cF$ be a class of graph polynomials.
A graph polynomial $F_0(G;\bar{x}) \in \cF$ is {\em universal for $\cF$} if every graph polynomial $F(G;\bar{x}) \in \cF$
can be written as a product of a $\VEC$-invariant and a substitution instance of $F_0(G;\bar{x})$.

A graph polynomial $F(G;\bar{x})$ is a {\em $\DC$-invariant (aka deletion-contraction invariant or chromatic invariant)}
if there are polynomials 
$g_1(\bar{x}), g_2(\bar{x}) \in \ZZ[\bar{x}]$
such that
$$
F(G;\bar{x}) =
g_1(\bar{x}) \cdot F(G_{-e};\bar{x}) + g_2(\bar{x}) \cdot F(G_{/e};\bar{x})
$$
and $F(G;\bar{x})$ is multiplicative, i.e.,
$$
F(G \sqcup G';\bar{x}) = F(G;\bar{x}) \cdot F(G';\bar{x}),
$$ 


\begin{theorem}
$Z(G;y,q) $
is universal for $\DC$-invariants.
\end{theorem}
For a proof see \cite{aigner2007course}.
There are quite a few similar theorems in the literature.

\subsection{Recipe theorems}
A {\em Tutte-Grothendieck invariant ($\TG$-invariant)} on graphs is a graph polynomial $F(G;\bar{x})$
such that
\begin{enumerate}[(i)]
\item
$F(G;\bar{x})$ is {\em multiplicative}
and
\item
there are polynomials 
$g_1(\bar{x}), \ldots,  g_5(\bar{x}) \in \ZZ[\bar{x}]$
and $F(G;\bar{x})$ satisfies the recurrence relation:
$$
F(G;\bar{x}) =
\begin{cases}
g_1(\bar{x}) & \text{ if } G = K_1 ;\\
g_2(\bar{x}) F(G-e) & \text{ if } e \text{ is a loop; } \\
g_3(\bar{x}) F(G/e) & \text{ if } e \text{ is a bridge; } \\
g_4(\bar{x}) \cdot F(G-e;\bar{x}) + g_5(\bar{x}) \cdot F(G_{/e};\bar{x}) & \text{ otherwise.}
\end{cases}
$$
\end{enumerate}
D. Welsh, \cite{welsh:94,welsh1999tutte}
coined the term {\em Recipe Theorem} for the following:
\begin{theorem}[T. Brylawski, 1972, \cite{brylawski1972decomposition}]
\label{thm:recipe}
The graph polynomial 
$$
U(G;w,q,v) = w^{m(G)} Z(G; q, \frac{q}{w})
$$ 
is a $\TG$-invariant and
any other TG-invariant is a substitution instance of $U(G;\bar{x})$ up to a prefactor.
\end{theorem}
In our terminology, this says that $U(G;\bar{x})$ is universal for $\TG$-invariants.
Theorem \ref{thm:recipe} is further discussed in
\cite{brylawski1992tutte,zaslavsky1992strong}.
The recipe here is the explicit definition of $U(G;w,q,v)$.

In his PhD thesis my former student, I. Averbouch, studied various other recurrence relations 
$\RR$ for polynomial graph invariants 
for which an $\RR$-universal graph polynomial exists. All his examples are multiplicative 
with respect to disjoint unions of graphs,
see \cite{phd:Averbouch}.
Among them we have
\begin{description}
\item[$U_{\DC}$] Generalizing the chromatic polynomial, based on edge deletion and edge contraction.
This is universal for chromatic ($\DC$)-invariants.
\item[$U_M$] Generalizing the matching polynomial, based on edge deletion and edge extraction.
\item[$U_{\EE}$] Generalizing both the matching and the chromatic polynomial, 
based on edge deletion, contraction and edge extraction,
This is universal for $\EE$-invariants, see
\cite{averbouch2010extension}.
\item[$U_{VE}$] Generalizing the Subgraph Component Polynomial from
\cite{tittmann2011enumeration} 
based in vertex deletion and vertex contraction and vertex extraction. 
\end{description}

The interlace polynomial was introduced in \cite{arratia2004interlace}.
Courcelle \cite{ar:Courcelle08} proved a Recipe Theorem generalizing the various interlace polynomials and the 
independence polynomial, based on
pivoting and local complementation.

\subsection{Recursive Definitions and Explicit Definitions of Graph Polynomials}

All the universal graph polynomials of the previous section not only have a recursive definition,
but also an explicit form as a graph polynomial.

Even if the class $\cF$ of graph polynomials has a recursive definition it may not have a universal
polynomial $U_{\cF}$. Nevertheless, under certain circumstances, the polynomials in $\cF$
may be of some specific form, say definable in second order logic $\SOL$.

\begin{tcolorbox}
Let $\cF$ be a class of graph polynomials which are well-defined by some (linear) recurrence relation. 
Then every $F \in \cF$ can be written explicitly as a polynomial
where the coefficients can be derived from the recurrence relation.
\end{tcolorbox}
Again, we consider such theorems as meta-theorems when the existence of the explicit 
polynomial and its coefficients
depend only on the form of the recursion in which $\cF$ is defined.

The technical apparatus needed to prove such a theorem is developed 
in \cite{godlin2012graph} and is rather involved. 
Due to lack of space, we refer the curious reader to \cite{godlin2012graph} in order to explore
how such a theorem can be formulated.

%% file: NMA-complexity.tex
\section{Complexity of Evaluating Graph Polynomials}
\label{se:complexity}
Another very early
meta-theorem in graph theory was the following observation
of Courcelle \cite{ar:CourcelleMosbah93,ar:Courcelle93}:
\begin{theorem}
Let $TW(k)$ be the class of graphs of tree-width at most $k$ and let $P$ be a graph property
definable in monadic second order logic $\MSOL$.
Then checking whether a graph $G \in TW(k)$ has property $P$ is in $\FPT$.
\end{theorem}
Here we think of this as a meta-theorem, because the conclusion only depends on the form of the formula
which defines $P$.


\subsection{Complexity: Fixed Parameter Tractability}

Our meta-theorem here is of the following form:
\begin{tcolorbox}
Let $GP(\cL)$ be the class of graph polynomials definable in the language of graphs in
some fragment $\cL$ of $\SOL$. 
Let $W(k)$ be the class of graphs of width at most $k$.
Then
evaluating 
$F(G; \bar{x}) \in GP(\cL)$ for $G \in W(k)$ is fixed parameter tractable
in the parameters $k$ and the $\cL$-formulas involved  in the definition of $F(G; \bar{x})$, irrespective
of the particular meaning of the formulas involved in defining $F(G; \bar{x})$.
\end{tcolorbox}
The width of a graph can be any notion of width for graphs as discussed in \cite{hlinveny2008width}.

The extension of Courcelle's Theorem to graph polynomials was part of the PhD theses of Udi Rotics \cite{phd:Rotics}
and later of Tomer Kotek \cite{phd:Kotek}, and published in the papers
\cite{courcelle2001fixed,ar:MakowskyTARSKI,ar:KotekMakowskyZilber08,kotek2013computational,kotek2018sequences}.

\begin{theorem}[Courcelle-type theorems for graph polynomials]
\label{th:courcelle-like}
Let $TW(k)$ be the graphs of tree-width at most $k$ and $GP(\MSOL)$ be the class of
$\MSOL$-definable graph polynomials. Let $F(G;\bar{x}) \in GP(\MSOL)$  and $\bar{a} \in \ZZ^s$.
Then evaluating $F(G;\bar{a})$ is in $\FPT$
in the parameters $k$ and the $\MSOL$-formulas involved  in the definition of $F(G; \bar{x})$.
\end{theorem}
Similar theorems are true for $CW(k)$, the graphs of clique-width at most $k$
and for graph polynomials in $GP(\cL)$, where $\cL$ ranges over logics which are variations of $\MSOL$.

\subsection{Complexity: Dichotomy Theorems}

A precursor to dichotomy theorems in complexity theory is Schaefer's theorem which shows that the boolean
satisfiability problem  is either $\bNP$-complete are in $\bP$ depending only on the syntactic description of
the class of formulas considered, \cite{ar:Schaefer78}. In our sense this can be also viewed as a meta-theorem.

The first dichotomy theorem for graph polynomials is due to Linial \cite{linial1986hard} and 
shows that evaluation of the chromatic polynomials $\chi(G;a)$ is 
$\sharp\bP$-hard for all complex values $a \in \CC$
with the exception of $a=0,1,2$.
This was generalized by Jaeger et al. in \cite{ar:JaegerVertiganWelsh90}.
It describes the complexity of evaluating the Tutte polynomial and shows that evaluation of $T(G; a, b)$
is $\sharp\bP$-hard for all complex values $(a,b) \in \CC^2$ 
with the exception of $(a,b) \in A$ where $A \subset \CC^2$ is an algebraic set of dimension $1$.

As a meta-theorem it takes the following form:
\begin{tcolorbox}
Let $\cF$ be a class of graph polynomials in $r$ indeterminates and $\cA_{\cF} \subseteq \CC^r$.
For every $F(G;\bar{x}) \in  \cF$ and every $\bar{a} \in \cA_{\cF}$ evaluating 
$F(G;a)$ is  $\sharp\bP$-hard. For $\bar{a} \not \in \cA_{\cF}$ it is in $\bP$.
\end{tcolorbox}
In
\cite{kotek2013computational}
we say that $F(G;\bar{x})$ has the {\em Difficult Point Property} if such an $\cA_{\cF}$ exists. 
See also \cite{goodall2018complexity}.

In the case of \cite{ar:JaegerVertiganWelsh90} $\cF$ is the class of $\TG$-invariants because the Tutte polynomial
is universal for $\TG$-invariants, and $\aA_{\cF}$
is an algebraic set in $\CC^r$ of dimension $< r$.
The same holds also for $\EE$-invariants and similar cases, \cite{phd:Hoffmann,phd:Averbouch}.

Let $A=A_H$ be the weighted adjacency matrix of a simple graph $H=(V(H),E(H))$ where $a_{i,j}$ is 
the weight of the edge $(i,j) \in E(H)$.
Given a graph $G = (V(G), E(G)$
we look at
$$
Z_A(G) = \sum_{\sigma:V \rightarrow [k]} \prod_{(i,j) \in E(G)} A_{\sigma(i), \sigma(j)}.
$$
If all the weights have value $1$, this counts the number of homomorphisms from $G$ to $H$.
If some of the values $a_{i,j}$ are treated as indeterminates we get graph polynomials in at most $n(G)^2$ variables.

In \cite{ar:BulatovGrohe2005} the following is shown:
\begin{theorem}[Bulatov and Grohe, 2005]
Let the weights of $H$ be non-negative  and real.
\begin{enumerate}[(i)]
\item
If $H$ is connected, then 
evaluating $Z_A(G)$ is in $\bP$ if $H$ is not bipartite and the row rank of $A_H$ is at most $1$,
or $H$ is bipartite and the row rank of $A_H$ is at most $2$. Otherwise it is $\sharp \bP$-hard.
\item
If $H$ is not connected, then evaluating $Z_A(G)$ is in $\bP$ if each of its
connected components satisfies the corresponding condition stated above.
Otherwise evaluating $Z_A(G)$ is $\sharp\bP$-hard.
\end{enumerate}
\end{theorem}
This generalizes \cite{ar:DyerGreenhill2000}.
There is a rich literature of similar dichotomy theorems in relation to Constraint Satisfaction Problems, 
see the references in \cite{ar:BulatovGrohe2005}. 

Further instances of the dichotomy meta-theorem can be found in 
\cite{blaser2007complexity,blaser2010complexity,goodall2018complexity}.

%% file: MA-unimodality.tex
\section{Graph Polynomials with Unimodal Coefficients}
\label{se:unimodality}
Read \cite{read1968introduction} conjectured already in 1968 that the absolute values of the coefficients of the
chromatic polynomial are unimodal. This has triggered interest in the question of the unimodality of the coefficents
of various graph polynomials. The conjecture was finally proved by Huh \cite{huh2012milnor}.
The case of the independence polynomial was discussed in \cite{alavi1987vertex}. More details follow in subsection
\ref{sse:bg}

The results in this section were inspired by the papers
\cite{barton2020acyclic,beaton2020unimodality,brown2018unimodality}, 
where the concept of almost unimodal graph polynomial was introduced, see subsection \ref{sse:au}.
These authors were discussing almost unimodality of particular graph polynomials 
such as the domination and independence polynomials, whereas we,
in V. Rakita's thesis \cite{phd:Rakita} and \cite{makowsky2023almost},
were interested in formulating meta-theorems. In this case, the meta-theorems have the following form:

\begin{tcolorbox}
Let $\fA$ be a family of graph properties given by a set of forbidden induced subgraphs $\Forb_{\fA}$
and $\cA \in \fA$.
Let $c^{\cA}_i(G)$ be the number of subsets $S\subseteq V(G)$ of a graph $G$ with $|S| =i$ such that
$G[S] \in \cA$.
Then the graph polynomials of the form
$$
F_{\cA}(G;x) = \sum_{i=0}^{d(G)} c^{\cA}_i(G) x^i
$$ 
share a certain behaviour which only depends on the existence of $\Forb_{\fA}$ 
irrespective of particular elements in $\Forb_{\fA}$.
\end{tcolorbox}

We first give some background on unimodality and almost-unimodality of polynomials
and only then state our meta-theorems.

\subsection{Real-rooted and unimodal graph polynomials}
\label{sse:bg}

Let $F(x) \in \R[x]$ be an univariate polynomial of degree $d$ with real coefficients, 
$$
F(x) = \sum_{i=0}^d a_i x^i.
$$
\begin{enumerate}[(i)]
\item
$F(x)$ is {\em real-rooted} if all its roots are in $\R$.
\item
The coefficients of
$F(x)$ are {\em log-concave} if for all $1 \leq j \leq d-1$
\\
{\em $a_j^2 \geq a_{j-1} a_{j+1}$.} 
\item
The coefficients of
$F(x)$ are {\em $\alpha$-concave} if for all $1 \leq j \leq d-1$
\\
{\em $a_j^2 \geq \alpha \cdot a_{j-1} a_{j+1}$.} 
\item
The coefficients of
$F(x)$ are {\em unimodal with mode $k$} if 
\\
{\em $a_i \leq a_j$ for $ 0\leq i < j \leq k$}
and
{\em $a_i \geq a_j$ for $ k \leq  i < j \leq d$.}
\item
$F(x)$ is {\em absolute unimodal with mode $k$ (log-concave)} if 
\\
the absolute values of $a_i$
are unimodal (log-concave).  
\item
These definitions, except for (i), apply to any sequence $a_i, 0 \leq i \leq d$,
even if it not interpreted as a sequence of coefficients of a polynomial.
\end{enumerate}

\begin{theorem}[Folklore]
\label{thm:NewtonsTheorem}
(i) implies (ii), (ii) implies (iv) and
none of the reverse implications holds.
\end{theorem}
The first part of the theorem is Newton's theorem. For a proof one may consult \cite{branden2015unimodality}.

Actually, Kurtz in \cite{kurtz1992sufficient} analyzed the situation further.
\begin{theorem}
\label{thm:Kurtz}
Let $F(x) \in \R[x]$ be an univariate polynomial of degree $d \geq 2$ with real coefficients,
$$
F(x) = \sum_{i=0}^d a_i x^i.
$$
\begin{enumerate}[(i)]
\item
If $F(x)$ is real-rooted, then its coefficients are $\alpha$-concave with
$$\alpha(d,i) = \frac{d-i+1}{d-i} \cdot \frac{i+1}{i}.$$
\item
If the coefficients of $F(x)$ are $4$-concave then all its roots are real and distinct.
\end{enumerate}
\end{theorem}

The sequence ${n \choose i}$ of the number of subsets of order $i$ of $V(G)$ is log-concave, hence unimodal.
More interestingly, 
let $m_i(G)$ be
the number of edge independent subsets (matchings) of $E(G)$ of order $i$.
The numbers $m_i(G)$ are also the coefficients of the matching generating polynomial
$$
M(G;x) =\sum_i m_i(G) x^i.
$$
\begin{theorem}
$M(G;x)$ is real-rooted, hence unimodal. 
\end{theorem}
There are two independent proofs of this.
It follows from the fact that all the roots of $M(G;x)$ are real for all graphs $G$, \cite{heilmann1970monomers},
see also \cite{gutman2016survey},
using Theorem \ref{thm:NewtonsTheorem}. Unimodality was also shown directly by A. Schwenk, \cite{schwenk1981unimodal}. 

Let $in_i(G)$, $0\leq i \leq n(G)$,
the sequence of the number of vertex independent subsets of $V(G)$ of order $i$. Denote by
$I(G;x) = \sum_i in_i(G)x^i$
the independence polynomial of $G$.
Real-rootedness, and unimodality of $I(G;x)$ has been studied extensively, see e.g. \cite{zhu2020unimodality}\cite{brown2018unimodality}\cite{cutler2017maximal} for some recent results and \cite{levit2005independence}\cite{levit2006independence} for a general introduction.

\begin{theorem}
\begin{enumerate}[(i)]
\item
$I(G;x)$ is not unimodal, 
\cite{alavi1987vertex}.
\item
For claw-free graphs the sequence of coefficients of $I(G;x)$ is real-rooted, hence unimodal, 
\cite{chudnovsky2007roots,bencs2014christoffel}.
\end{enumerate}
\end{theorem}
However, it is easily seen that the set of counterexamples $G$ for (i)
given in  \cite{alavi1987vertex}
has measure $0$ among the random graphs $\mathcal{G}(n,p)$.
Similarly for (ii),
it is easily seen that the claw-free graphs have measure $0$ among the random graphs $G(n,p)$.

This leaves open whether $I(G;x)$ is real-rooted, or at least unimodal, for other graph classes. 
Specifically, we may ask whether the {\em independence polynomial $I(G;x)$} is unimodal for ``most graphs'', 
in the following sense:

\begin{definition}
Let $P$ be a graph polynomial. We say $P$ is {\em almost unimodal} if for almost all graphs $G$ the polynomial 
$P(G;x)$ is unimodal. In other words, $P$ is almost unimodal if 
for random graphs $G \in \mathcal{G}(n, 1/2)$ we have
$$
\lim_{n\goesto \infty} \PP(P(G(n,1/2),x) \text{ is unimodal} )=1.
$$
\end{definition}

\begin{problem}
Is $I(G;x)$  almost unimodal?
\end{problem}

Let $\chi(G,x) = \sum_i c_i(G) x^i$ be the chromatic polynomial of $G$.
The case of the chromatic polynomial of a graph is slightly different.
The sequence $c_i(G)$ is alternatingly positive and negative. However,
it was conjectured by R.C. Read, \cite{read1968introduction}, that the absolute values 
$|c_i(G)|$ form a unimodal sequence.
J. Huh, \cite{huh2012milnor} finally proved the conjecture.
\begin{theorem}[J. Huh, 2012]
For every graph $G$ the chromatic polynomial $\chi(G,x)$ is absolute unimodal.
In fact the sequence $|c_i(G)|$ is log-concave.
\end{theorem}
\subsection{Almost Unimodality of the Coefficients}
\label{sse:au}
Let 
$$
F(G;x) = \sum_{i=0}^{d(G)} a_i(G) x^i
$$ 
be a graph polynomial.
$F(G;x)$ is (almost) unimodal if for (almost) all graphs $G$, the sequence $a_i(G)$ is unimodal.

Recall that for a graph property $\cA$
$c^{\cA}_i(G)$ is the number of subsets $S\subseteq V(G)$ of a graph $G$ with $|S| =i$ such that
$G[S] \in \cA$.

The following two theorems are proved in
\cite{makowsky2023almost}.
\begin{theorem}[Almost Unimodality Theorem]
\label{th:unimodal}
If $\cA$ is a non-trivial co-hereditary graph property then
for almost all graphs $G$, the sequence $c^{\cA}_i(G), 0\leq i \leq n(G)$ is unimodal.
In other words
$$
F_{\cA}(G;x) = \sum_{i=0}^{d(G)} c^{\cA}_i(G) x^i
$$ 
is almost unimodal.
\end{theorem}

$F(G;x) = \sum_{i=0}^{i=d(G)} a_i(G) x^i$ is real-rooted if all its roots are real.
Real-rooted univariate polynomials are unimodal.

\begin{theorem}[Real-rootedness Theorem]
\label{th:realrooted}
Let $\cA$ be hereditary and such that it contains a graph $G \in \cA$ which is neither a clique nor an edgeless graph.
Then $F_{\cA}(G;x)$ is real-rooted if and only if $G \in \cA$.
\end{theorem}

In \cite[Theorem 4.17]{makowsky2014location}
it is shown that for every univariate graph polynomial $F(G;x)$ with integer coefficient
there is a 
graph polynomial $F_1(G;x)$
with integer coefficients 
with the same distinguishing power\footnote{
For a precise notion one can use
s.d.p.-equivalence,
see Section \ref{se:semantic}.
}
such that all the roots of  $F_1(G;x)$ are real.

%% file: MA-harary.tex
\section{Harary polynomials}
\label{se:harary}

A graph property $\cA$ is {\em trivial} if it is empty, finite, cofinite, 
or it contains
all finite graphs
(up to isomorphisms). 
Let $\cA$ be a non-trivial graph property.

Harary  polynomials are of the form
\begin{gather}
\chi_{\cA}(G; x) = \sum_{i \geq 1} b_i^{\cA}(G) x_{(i)},
\end{gather}
where $b_i^{\cA}(G)$ is the number of partitions of $V(G)$ into $i$ non-empty parts, where each part induces a graph in ${\cA}$,
and $x_{(i)}$ is the falling factorial.

Harary polynomials are special cases of generalized chromatic polynomials as described in
\cite{ar:MakowskyZilber2006,ar:KotekMakowskyZilber08,ar:KotekMakowskyZilber11}.
They are the counting analogues of conditional colourings as defined in \cite{harary1985conditional}.

A meta-theorem could be formulated as follows:
\begin{tcolorbox}
Let $\frak{P}$ a class of graph properties such as hereditary, monotone or minor-closed classes.
If $\cA \in \frak{P}$ then the Harary polynomials $\chi_{\cA}(G; x)$
exhibit a certain behaviour.
\end{tcolorbox}


The main questions in  V. Rakita's thesis and in \cite{phd:Rakita,herscovici2020harary} 
ask whether Courcelle's Theorem and its variations
can be applied to Harary polynomials for non-trivial graph properties. 
This amounts to asking: 
\begin{enumerate}[(i)]
\item
Are there non-trivial graph properties $\cA$ such that the Harary polynomial $\chi_{\cA}(G,x)$ is  
$\MSOL_g$-definable? 
\item
Are there non-trivial graph properties $\cA$ such that the Harary polynomial $\chi_{\cA}(G,x)$
is an $\EE$-invariant and hence $\MSOL_h$-definable?
\end{enumerate}

Recall that a graph property $\cA$ is {\em hereditary (monotone, minor-closed)}
if it is closed under taking induced subgraphs (subgraphs, minors).
Clearly, 
if $\cA$ is minor-closed, it is also monotone, and
if $\cA$ is monotone, it is also hereditary.

A graph property $\cA$ is {\em ultimately clique-free} if there exist $t \in \NN$ such that no graph $G \in \cA$
contains a $K_t$, i.e., a complete graph of order $t$.
Analogously, $\cA$ is {\em ultimately  biclique-free} if there exist $t \in \NN$ such that no graph $G \in \cA$
has $K_{t,t}$  as a subgraph (not necessarily induced). $K_{t,t}$  is the  complete bipartite graph on two sets of size $t$.
Clearly, ultimately biclique-free implies ultimately clique-free, but not conversely.

\begin{theorem}
\label{th:main-1}
Let  $\cA$ be a graph property and $\chi_{\cA}(G,x)$ the Harary polynomial associated with  $\cA$.
\begin{enumerate}[(i)]
\item
If $\cA$ is hereditary, monotone, or minor closed,
then
$\chi_{\cA}(G;x)$ is an $\EE$-invariant if and only if 
$\chi_{\cA}(G;x)$ 
is the chromatic polynomial $\chi(G;x)$.
\item
If $\cA$ is ultimately clique-free (biclique-free),
$\chi_{\cA}(G;x)$ is not $\MSOL_g$-definable.
\end{enumerate}
\end{theorem}


%% file: MA-weakly.tex
\section{Weakly Distinguishing Polynomials}
\label{se:weakly}
This section is based on V. Rakita's MSc thesis, \cite{msc:Rakita} and on \cite{makowsky2019weakly}.
It was inspired by the Conjecture of 
Bollob\'as, Pebody and Riordan (Conjecture \ref{BPR}) below, in attempt to identify graph polynomials $P$
where there are many  graphs $G$ which are not $P$-unique.
We follow the terminology of \cite{makowsky2019weakly}.

\subsection{Weakly distinguishing on a graph property $\cC$}

For a graph property $\mathcal{C}$, denote by $\mathcal{C}(n)$ the graphs of order $n$ in $\mathcal{C}$. 
We only consider properties such that $\mathcal{C}(n)$ is non-empty for all sufficiently large $n$.

Let $P$ be a graph polynomial. We say that two non-isomorphic graphs $G$ and $H$ are 
$P$-mates if $P(G)=P(H)$, and that $G$ is $P$-unique if it has no $P$ mates. 
$P$ is {\em trivial} if all graphs $G,H$ are $P$-mates.
$P$ is {\em complete} if all graphs $G$ are $P$-unique.

In this section we investigate conditions which imply that almost all graphs in 
$\mathcal{C}$ have a $P$-mate. More formally, we give the following definitions:

Let $\mathcal{G}(n)$ be the family of graphs 
of order $n$ with $V(G)=\{1,...n\}$. 
Let $P$ be a graph polynomial and 
let $U_P(n)$ be the set of 
$P$-unique graphs in $\mathcal{G}(n)$.
\begin{definition}
$P$ is {\em weakly distinguishing} if 
$$\lim_{n\goesto \infty} \dfrac{|U_P(n)|}{|\mathcal{G}(n)|}=0$$
and that $P$ is {\em almost complete} if 
$$\lim_{n\goesto \infty} \dfrac{|U_P(n)|}{|\mathcal{G}(n)|}=1$$
\end{definition}

In \cite{bollobas2000contraction},
Bollob\'as, Pebody and Riordan conjectured:

\begin{conjecture}[BPR-conjecture]
\label{BPR}
The  chromatic and Tutte polynomials are almost complete. 
\end{conjecture}
In \cite{makowsky2019p}, the analogous question 
for $r$-regular hypergraphs 
was considered, and for $r \geq 3$ the conjecture was refuted.
In \cite{noy2003graphs} it was observed, as a remark in the conclusions, that
the independence polynomial $In(G;x)$ 
is weakly distinguishing on all finite graphs.
Clearly, we are intersted in identifying $\cC$ as large as possible.

The meta-theorems I have in mind can be of two forms:

\begin{tcolorbox}
Let $\cF$ be a class of graph polynomials and $\cA_{\cF}$ be a family of graphs.
Then almost every $G \in \cA_{\cF}$ is $F$-unique.
\end{tcolorbox}

\begin{tcolorbox}
Let $\cF$ be a class of graph polynomials and $\cA_{\cF}$ be a family of graphs.
Then almost every $G \in \cA_{\cF}$ has a $F$-mate.
\end{tcolorbox}
To make these into meta-theorems, $\cF$ or $\cA_{\cF}$  should be given in a 
framework of logic or of recursive definitions.

\begin{problem}
Is the universal $\EE$-invariant $U_{\EE}$ from section \ref{se:recursion} almost complete?
\end{problem}

\subsection{Small Addable Classes of Graphs}
In this subsection the terminology is from \cite{mcdiarmid2005random,mcdiarmid2006random}.

\begin{definition}
\begin{enumerate}[(i)]
\item
We say a graph property $\mathcal{A}$ is {\em decomposable} if it is closed under disjoint union, 
and for all $G\in \mathcal{A}$ every component of $G$ is in $\mathcal{A}$.
\item
We say a graph property $\mathcal{A}$ is {\em bridge addable} if for each graph $G\in \mathcal{A}$ 
and every two vertices $u,v$ in different components of $G$ the graph obtained 
by adding an edge between $u$ and $v$ is also in $\mathcal{A}$.
\item
We say a graph property $\mathcal{A}$ is {\em addable} if it is decomposable and bridge addable.
\item
Let $\aA$ be a graph property, and denote by $\aA_n$ the graphs of order $n$ in $\aA$. 
We say that a graph property $\aA$ is {\em small} if there exists a constant $a>0$ 
such that $|\aA_n|\leq a^nn!$ for all sufficiently large $n$.
\end{enumerate}
\end{definition}
Being small is a purely technical condition needed in \cite{mcdiarmid2006random}.
In \cite{norine2006proper} it is shown that a proper minor closed graph property is small.

In \cite{rakita2019weakly} it was proven that an infinite number of graph polynomials, 
among them the independence, clique and harmonious polynomials, are weakly distinguishing on the class
of small addable graphs.

In particular the following theorem is shown:

\begin{theorem}
Let $F(G; \bar{x})$ be  a chromatic invariant, a $\TG$-invariant, or even an $\DCE$-invariant. 
Then $F(G; \bar{x})$ is weakly distinguishing on the class
of small addable graphs.  In particular this holds also for the class of
proper minor-closed addable graphs. 
\end{theorem}

\ifskip\else
\subsection{Rakita's suggestion}

Correction: 
$\cC$ needs to be that every connected graph in $\cC$ is is $k$ connected, and $\cC$ is closed under dusjoint unions.

Let $\cC$ be a graph property such that all graphs in $\cC$ are $k$ connected for $k>1$, and let $H_{\cC}(G;x)$ be its Harary polynomial. 
Proposition: 
Let $G$ be a graph, and let $A,B$ be a partition of the vertices of $G$ such that $G[A]$ and $G[B]$ are $k$ connected for $k>1$ , and there is at most one edge between $A$ and $B$. Then 
$H_{\cC} (G;x)= H_{\cC} (G[A];x) H_{\cC} (G[B];x)$
Corollary: Let $\cA$ be a non empty small addable property, and $\cC$ and $H_{\cC}(G;x)$  as above. Then $H_{\cC}(G;x)$  is weakly distinguishing on $\cA$.

\fi 

\subsection{Weakly distinguishing Harary polynomials}

\begin{definition}[R-class]
Let $\cC$ be a graph property. 
$\cC$ is an $R$-class if
\begin{enumerate}[(i)]
\item
$\cC$ is closed under disjoint unions, and
\item
every connected component of $G \in \cC$ is at least $2$-connected.
\end{enumerate}
\end{definition}

\begin{examples}
The following are $R$-classes of graphs:
\begin{enumerate}[(i)]
\item
The closure $DU(\cC)$ under disjoint unions of any class $\cC$ where all its connected members are at least $2$-connected.
In particular, when $\cC$ consists of all complete graphs $K_n, n \geq 3$, grids $Grid_{m,n}: m, n \geq 2$, 
wheels $W_n: n \geq 2$, etc.
\item
Any addable class where all its connected members are at least $2$-connected.
\end{enumerate}
\end{examples}

\begin{definition}[R-Harary polynomial]
The Harary polynomial $H_{\cC}(G;x)$  is an $R$-Harary polynomial if $\cC$ is an $R$-class.
\end{definition}

\begin{proposition}
Let  $H_{\cC}(G;x)$  be an $R$-Harary polynomial. 
Let $G$ be a graph, and let $A,B$ be a partition of the vertices of $G$ such that $G[A]$ and $G[B]$ 
are at least $2$-connected,  and there is at most one edge between $A$ and $B$. 
Then 
$$
H_{\cC} (G;x)= H_{\cC} (G[A];x) \cdot H_{\cC} (G[B];x).
$$
\end{proposition}

\begin{corollary} 
Let $\cA$ be a non-empty small addable property, and $\cC$ and $H_{\cC}(G;x)$  as above. 
Then $H_{\cC}(G;x)$  is weakly distinguishing on $\cA$.
\end{corollary}

%% file: NMA-semantic.tex
\section{Semantic vs syntactic properties of graph polynomials}
\label{se:semantic}

Let $\cA$ be a graph property.
Two graph polynomials $F_1(G, \bar{x}), F_2(G, \bar{x})$ are {\em $(s.d.p)$-equivalent} om $\cA$, 
written as $F_1 \sim_{s.d.p}^{\cA} F_2$,
if for all similar graphs $G_1, G_2 \in \cA$ we have
$$
F_1(G_1, \bar{x}) = F_(G_2, \bar{x})
\text{ if and only if }
F_2(G_1, \bar{x}) = F_2(G_2, \bar{x}).
$$
$(s.d.p)$-equivalent stands for Distinguishing Power on Similar graphs.

\begin{example}
Let $m_k(G)$ be the number of  $k$-matchings of a graph $G$.
The matching generating polynomial
$$
m(G;x) =  \sum_k m_k(G) x^k
$$
and the matching defect polynomial 
$$
M(G;x) =  \sum_k (-1)^k m_k(G) x^{n-2k}
$$
are for every graph $G$ different polynomials, but they are $(s.d.p)$-equivalent.
\end{example}

Let $\cF$ be a family of graph polynomials and $\hat{\cF}$ be its closure under
$(s.d.p)$-equivalence, i.e., $\cF \subseteq \hat{\cF}$, and if  $F \in \hat{\cF}$ and $F'$ is $(s.d.p)$-equivalent to $F$
then $F' \in \hat{\cF}$.
We call a true statement $\Phi$ about $\cF$ a {\em semantic statement} if $\Phi$ also holds for $\hat{\cF}$.

\begin{example}
Let $\cF$ be a family of graphs and $\cA$ be a graph property.

We look at the following statements:
\begin{enumerate}[(i)]
\item
For every $F_1, F_2 \in \cF$ and every $G \in \cA$
we have $F_1(G;\bar{x}) = F_2(G;\bar{x})$.
\item
For every $F_1, F_2 \in \cF$ 
we have $F_1(G;\bar{x}) \sim_{s.d.p}^{\cA} F_2(G;\bar{x})$.
\end{enumerate}
(i) is not a semantic property of $\cF$, whereas (ii) is.

Both statements are true when $\cF$ just contains the characteristic polynomial of the adjacency matrix,
and the matching defect polynomial
and $\cA$ are the forests, \cite{godsil1981theory} and \cite[Chapter 2, Corollary 1.4]{bk:Godsil93}.
\end{example}

The distinction between semantic and syntactic statements about graph polynomials was  developed in the author's
graduate courses on graph polynomials starting in 2005 and first explicitly discussed and published
in \cite{makowsky2014location}. There it is shown that statements about the location of the roots of graph polynomials
are syntactic, see also \cite{makowsky2019logician}. Using the same methods one can show that
statements about universality, complexity, unimodality or being a Harary polynomial are syntactic. 
In contrast to this, statements 
about completeness, almost-completeness, weakly distinguishing, etc are semantic statements.

%% file: NMA-conclu.tex
\section{Conclusions and Suggestion for Further Research}
\label{se:conclu}

We have discussed graph polynomials from the point of view of meta-theorems.
Peter Tittmann's encyclopedic book \cite{bk:Tittmann} surveys the landscape of graph polynomials.
The recent paper by Graham Farr and Kerri Morgan, \cite{ar:FarrMorgan2024}, somehow bridges the two points of view.

We now discuss possible directions of further research towards more meta-theorems.

\paragraph{\bf Universal polynomials and recurrence relations}
\begin{problem}
Characterize recurrence relations defining a class of graph polynomials which have a universal polynomial.
\end{problem}

\begin{problem}
Characterize recurrence relations defining a class of graph polynomials such that each of its graph polynomials
can be written in an explicit  form. In other words, simplify the framework of \cite{godlin2012graph}  
\end{problem}

\paragraph{\bf Complexity of evaluation}
\begin{problem}
Are there other graph parameters besides variants of width for which Courcelle-like
theorems like Theorem \ref{th:courcelle-like} can be established?
\end{problem}

\begin{problem}
Characterize classes of graph polynomials which satisfy a dichotomy theorem.
\end{problem}

\paragraph{\bf Unimodality}
If $A \subseteq V(G)$ is independent in $G$, so are the subsets of $A$.
Let $\cA \subseteq \wp(V(G)$ be closed under subsets.
Let $c^{\cA}_i(G)$  be the number of sets of size $i$ in $\cA$.
We define a graph polynomial 
$$
P^{\cA}(G;x) =  \sum_i c^{\cA}_i(G) x^i.
$$
\begin{problem}
Under what condition on $\cA$, besides being closed under subsets, is $P^{\cA}(G;x)$ almost unimodal?
In particular, is the Independence polynomial almost unimodal.
\end{problem}

\paragraph{\bf Harary polynomials}
Let $\fA$  a family of graph properties and let for $\cA \in \fA$
$\chi_{\cA}(G; x) = \sum_{i \geq 1} b_i^{\cA}(G) x_{(i)}$  be its Harary polynomials.
\begin{problem}
Find conditions for large families $\fA$ such that all its Harary polynomials 
are mutually incomparable in distinguishing power.
\end{problem}

\paragraph{\bf Weakly distinguishing graph polynomials}
\begin{problem}
Is the most general $\EE$-invariant almost complete?
\end{problem}

\begin{problem}
Is there a $\SOL$-definable ($\MSOL$-definable) graph polynomial which is almost complete?
\end{problem}

\paragraph{\bf Semantic vs syntactic statements}

We have noted that most of our meta-theorems are of syntactic nature, i.e., they make statements which refer to the syntactic form
in which the graph polynomials are represented, rather than to their distinguishing power.

\begin{problem}
Find more meta-theorems which are semantic statements.
\end{problem}